\documentclass[12pt, oneside,reqno]{amsproc}   	% use "amsart" instead of "article" for AMSLaTeX format
\usepackage{geometry}                		% See geometry.pdf to learn the layout options. There are lots.
\usepackage{graphicx}				% Use pdf, png, jpg, or eps with pdflatex; use eps in DVI mode
\usepackage{amssymb}
\usepackage{amsmath}
\usepackage[all]{xy}
\usepackage{mathtools}
\usepackage{tikz-cd}
\usepackage{enumitem}

\usepackage{tikz-cd}

\newtheorem{thm}{Theorem}[section]
\newtheorem{lem}[thm]{Lemma}
\newtheorem{cor}[thm]{Corollary}

\newtheorem{conj}{Conjecture}[section]
\newtheorem{claim}{Claim}[section]

\numberwithin{equation}{section}

\def\Pb{\ifmmode{\Bbb P}\else{$\Bbb P$}\fi}
\def\Z{\ifmmode{\Bbb Z}\else{$\Bbb Z$}\fi}
\def\C{\ifmmode{\Bbb C}\else{$\Bbb C$}\fi}
\def\R{\ifmmode{\Bbb R}\else{$\Bbb R$}\fi}
\def\S{\ifmmode{S^2}\else{$S^2$}\fi}

\def\S{\cal S}

\newenvironment{pf}{\paragraph{Proof:}}{\hfill$\square$ \newline}

\begin {document}
	
\title{On the topological rigidity of compact self shrinkers in $\R^3$}
\begin{abstract}In this note we show that self shrinkers in $\R^3$ are ``topologically standard" in that any genus $g$ compact self shrinker is ambiently isotopic to the standard genus $g$ embedded surface in $\R^3$. As a consequence self shrinking tori are unknotted.  \end{abstract}

\author {Alexander Mramor \& Shengwen Wang}
\address{Department of Mathematics, University of California Irvine, Irvine, CA 92617}
\address{Department of Mathematics, Johns Hopkins University, Baltimore, MD 21218}
\email{mramora@uci.edu,swang@math.jhu.edu}

\maketitle

\section{Introduction.}

Self shrinkers, which in this article refers to hypersurfaces $M^n \subset \R^{n+1}$ satisfying
\begin{equation}
H + \frac{1}{2}x^\perp = 0
\end{equation} 
are singularity models for the mean curvature flow and hence are of considerable interest to the study of the flow. The mean convex case, $H \geq 0$, is understood very well; by Huisken \cite{H} all mean convex self shrinkers are round generalized cylinders $S^{n-r}(\sqrt{2(n-r)}) \times \R^r$ for some $0 \leq r \leq n$ (indeed, the only mean convex compact self shrinker is the round sphere and namely has genus 0). For $n =2$ Brendle showed the only self shrinker of genus 0, mean convex or not, was the round sphere of radius 2. The general case is much less well understood however, and only a patchwork of non mean convex embedded self shrinkers such as Angenent's Torus \cite{Ang}, gluing constructions by Kapoluleas, Kleene, and Moller \cite{KKM}, minimax surfaces by Ketover \cite{Ket}, and numerical evidence by Chopp in \cite{C} are known. For a more complete overview of the landscape of self-shrinkers, see the excellent recent survey by Drugan, Lee, and Nguyen \cite{DLM}.
$\medskip$

Of course conceivably there are many more examples of self shrinkers exist than those described above. In this note however we show that at any rate when $n=2$ they must be "topologically standard" in the following way: 
\begin{thm}\label{main} Let $F: \Sigma_g \to \R^3$ be an embedded closed self shrinker of genus $g$. Then it is isotopic to the standard genus $g$ surface in $\R^3$. 
\end{thm} 

As a corollary:
\begin{cor} Let $F: T^2 \to \R^3$ be an embedded self shrinking torus. Then it is unknotted. 
\end{cor}

We recall a knot is a closed embedded curve $\gamma: S^1 \to \R^3$. $\gamma$ is unknotted if it is ambiently isotopic to the equator of the round 2-sphere of radius 1 in $\R^3$. Similarly, in this article we will say a hypersurface $F: T^2 \to \R^3$ is $\textbf{unknotted}$ if it is ambiently isotopic to a tubular neighborhood of an unknotted curve $\gamma$. 
$\bigskip$

$\textbf{Acknowledgements.}$ The authors thank their advisors, Richard Schoen and Jacob Bernstein respectively, for their guidance. The authors also thank the Fields institute, particularly for hosting the July 2017 geometric analysis summer school, where the authors first met and first discussed the topic of this paper. Finally the authors thank Hojoo Lee for discussing with us the history of the Angenent torus uniqueness conjecture. 

\section{Proof of Theorem \ref{main}.}

Let's start by setting some notation: $M$ will be a compact self shrinker of genus $g$. Denote by $N_{int}$ and $N_{ext}$ the interior and exterior components of $M$ in $\R^3$, and $\widetilde{N_{int}}$ and $\widetilde{N_{ext}}$ by their corresponding sets in $S^3$ by adding $\{\infty\}$ to $\R^3$ and thinking of $S^3$ as the one point compactification of $\R^3$ via stereographic projection (since $M$ is compact, as sets note $\widetilde{N}_{int} = N_{int}$).  
$\medskip$

We will also need that the self-shrinker $M$ is a minimal surface in the Gaussian metric $(\mathbb R^3,\frac{1}{(4\pi)}e^{\frac{-|x|^2}{4}}\delta_{ij})$,  so as  $M$ is the boundary of both $N_{int}$ and $N_{ext}$ they both are manifolds with mean convex boundary.
$\medskip$

In \cite{L}, Lawson proved the analogue of theorem 1 for minimal surfaces in $S^3$ by showing that minimal surfaces are Heegaard splitings of $S^3$; which is to say their interior and exterior components are Handlebodies. He then appealed to a deep theorem of Waldhausen \cite{W} that says any two genus $g$ Heegard spliting of the 3 sphere are ambiently isotopic. For the main geometric step in his argument, which is essentially lemma 2.1 below, Lawson uses some facts about the round metric which the Gaussian metric doesn't satisfy (we discuss this more in the concluding remarks). Instead we apply a minimization argument of Frohman-Meeks \cite{FM}; some modifications similar to what's found in Brendle's paper \cite{B1} are necessary though since the metric decays to the zero form at infinity. 
$\medskip$

\begin{lem}\label{surjectivity} Let $M$ be a connected compact self shrinker in $\R^3$ or equivalently a connected minimal surface in $(\R^3,G)$ where $G$ is the Gaussian metric. Then the map
\begin{equation}
 \pi_1(\partial N_{ext}) \xrightarrow{i_*} \pi_1(N_{ext})
 \end{equation}
 induced by inclusion is surjective. The same is true for $N_{int}$. 
 \end{lem}
 
In the following we state precisely the minimization lemma of Frohman and Meeks we mentioned above; in the statement an almost-complete Riemannian surface means namely a complete metric space with respect to the distance function induced by infimum of length curves joining 2 points. It's proof is sketched below in the course of proving surjectivity $\pi_1(\partial N_{ext}) \xrightarrow{i_*} \pi_1(N_{ext})$:
 \begin{lem}\label{FMlemma}(Lemma 3.1 of \cite{FM})
 
Suppose $N$ is a connected, orientable, almost-complete Riemannian three-manifold with more than one boundary component. If $\partial N$ has nonnegative mean curvature with respect to the outward pointing normal, then $N$ contains a properly embedded, orientable, least-area minimal surface.
 \end{lem}
 $\medskip$
 
 Before starting the proof, let's recall some basic facts about the length functional. For a path $\gamma$ in a Riemannian manifold $(N,\mu)$, denote its length with respect to $\mu$ by $L_\mu(\gamma) = \int \mu(\dot{\gamma}, \dot{\gamma}) dt$. This let's us induce a metric space topology on $(N, \mu)$ by defining the distance between two points $p$ and $q$ to be given by $\inf\limits_{\{ \gamma \mid \gamma(0) = p, \gamma(1) = q\}}L_{\sigma \mu}(\gamma)$. We note that if, as positive definite matrices over each point, two metrics $g_1, g_2$ on $N$ satisfy $cg_1 < g_2 < Cg_1$ for $c, C > 0$, then we have for their length functionals $L_{cg_1} < L_{g_2} < L_{Cg_1}$. Hence if as a metric space induced by the length functional $(N, g_1)$ is complete, so is $(N, g_2)$. 
 $\medskip$

\begin{pf} Since $M$ is compact and in every compact ball there is some constant $\sigma$ so that $\sigma \delta_{ij} < G < \delta_{ij}$, from the comments above the $(N_{int}, G)$ is almost complete so we get from lemma \ref{FMlemma} a stable self shrinker $\pi(\Sigma) \subset N_{int}$. From here we may proceed from claim 2.1 below onward to conclude. 
$\medskip$

For $N_{ext}$ we follow closely the argument of Frohman and Meeks (but enough modifications are necessary that we don't just quote their lemma). Suppose $\pi_1(M)\rightarrow\pi_1(N_{ext})$ is not surjective, $M=\partial N_{ext}$. Then by elementary covering space theory $\partial \bar N_{ext}$ is not connected, where $\partial\bar N_{ext}$ is the boundary of the universal cover $\pi:\bar N_{ext}\rightarrow N_{ext}$. Denote by $\partial_1,\partial_2$ two (of perhaps many) connected components of $\partial \bar N_{ext}$. 
$\medskip$

We know $\partial N_{ext}=M$ is a minimal surface in $(N_{ext},G=e^{-\frac{|x|^2}{4}}\delta_{ij})$, so $\partial_1,\partial_2$ as subsets of the lift of $M$ are also minimal surfaces in $(\bar N_{ext}, \bar G)$, where $\bar G$ is the lift of $G$ and hence are mean convex.  Denote $\gamma\subset\bar N_{ext}$ a smooth curve connecting $\partial_1,\partial_2$. At this point Frohman and Meeks find a stable minimal surface disjoint from the boundary but $(\bar N_{ext}, \bar G)$ is not almost-complete so we can't do that yet, so we perturb the metric as in \cite{B1}. Let $\Psi$ be a smooth bump function in $\R^3$ defined by
\begin{equation}
\begin{split}
\Psi(x)&=1,x\in B(0,1)\\
\Psi(x)&=0,x\in \mathbb R^3\setminus B(0,2)\\
\end{split}
\end{equation}
and
\begin{equation}
\Psi_k(x)=1-\Psi(\frac{x}{k})
\end{equation}
where $B(a,R)$ is the open ball of radius $R$ centered at $a$. We use these functions to perturb the Gaussian metric $G$ of $N_{ext}$ and define for any $k\in\mathbb N$
\begin{equation}
G_k=(e^{-\frac{|x|^2}{4}}+\Psi_k)\delta_{ij}
\end{equation}

We can choose $k_0\in\mathbb N$ large enough so that
\begin{equation}
\begin{split}
&M\subset B(0,k_0)\\
&\pi(\gamma)\subset B(0,k_0)
\end{split}
\end{equation}

Now for any $k>k_0$, $\partial_1,\partial_2$ are minimal surface in the lifted metric $(\bar N_{ext},\bar G_k)$, moreover $(\bar N_{ext},\bar G_k)$ is a complete Riemannian metric and an almost complete metric space (from what we discussed in the paragraph before the start of the proof) satisfying conditions of Lemma \ref{FMlemma}, so we get a properly embedded connected area minimizing surface $\Sigma_k\subset(\bar N_{ext},\bar G_k)$. For the sake of completeness let's briefly explain how $\Sigma_k$ is found:
$\medskip$

Let $\Pi_1 \subset \Pi_2 \subset \ldots$ be a smooth compact exhaustion of $\partial_1$ with $p \in \Pi_1$. Recalling that the metric was perturbed away from $M$ so the boundary in the universal cover is still mean convex, by standard existence and regularity theory for minimizers (completeness of space necessary) for each $i$ we may find $\Sigma_i$ properly embedded, orientable area minimizer with boundary $\partial \Pi_i$.
$\medskip$

 Since each of the $\Sigma_i$ are area minimizing we get local area bounds which lets us apply standard compactness and regularity theory to take a limit of the $\Sigma_i$ to get an area minimizing surface we denote $\Sigma_k$. Since each of the $\Sigma_i$ intersected the arc $\gamma$ with odd intersection number, the limit surface $\Sigma_k$ is nonempty and intersects with $\gamma$. To proceed we need the following sub-lemma:

\begin{claim}
$\Sigma_k$ is disjoint from $\partial\bar N_{ext}$
\end{claim}
\begin{pf}
If not, then by maximum principle $\Sigma_k\subset\partial\bar N_{ext}$ and $\Sigma_k$ agrees with one of the connected components of $\partial\bar N_{ext}$, thus $\pi(\Sigma_k)=M$.
$\medskip$

Since $\Sigma_k$ is area minimizing in $(\bar N_{ext,},\bar G_k)$, $M$ is stable in $(N_{ext}, G_k)$. But $k>k_0$ implies $G$ and $G_k$ agree in a neighborhood of $M$, meaning $M$ is stable in $(N_{ext},G)$. Now recall from Colding-Minicozzi (namely theorem 4.14 in \cite{CM}) that if we take $L_M=\Delta_M+|A|^2+\frac{1}{2}-\frac{1}{2}\nabla_{x^T}$, then for a self shrinker the second variation is given by\begin{equation} 
- \int_M e^{-\frac{|x|^2}{4}} f Lf = \int_M e^{- \frac{|x|^2}{4}}( | \nabla^M f|^2) - |A|^2 f^2 - \frac{1}{2} f^2)
\end{equation}
for test functions that vanish along the boundary of $M$. If $M$ is closed and we plug in the test function $f=1$ we have $L(f)>0$ so the constant outward normal vector field on $M$ will decrease area, implying that $M$ is not stable in $(N_{ext},G)$, a contradiction. Similarly for $N_{int}$ we may plug in $f = -1$ to see the $\Sigma$ we find from just applying lemma 2.2 directly (see start of proof) is disjoint from $\partial \overline{N}_{int}$.

\end{pf}

Now $\Sigma_k\cap B(0,k_0)\neq\emptyset$ because they all intersect with $\gamma$. We now let let $k\rightarrow\infty$. The ambient manifolds $(\bar N_{ext}, \bar G_k)\rightarrow (\bar N_{ext},\bar G)$ and the stable minimal surfaces after passing to a subsequence (uniform ambient curvature bounds imply uniform local area bounds for area minimizers, which implies curvature (and higher order) bounds) converge to some surface $\Sigma_k \to \Sigma$ where $\Sigma\subset(\bar N_{ext},\bar G)$ is a non-empty stable minimal surface that is disjoint from $\partial\bar N_{ext}$, provided $\partial\bar N_{ext}$ is disconnected. When projecting down to the base, $\pi(\Sigma)\subset(N_{ext},G)$ is an embedded minimal surface that is disjoint from $M$, namely it's a self-shrinker in $\mathbb R^3$ and lies in $N_{ext}$.
$\medskip$ 

As disjoint self-shrinkers, we have the distance between $M$ and $\pi(\Sigma)$ becomes 0 at $t = 0$ (equation (1.1) is a normalized equation for self shrinkers that implies solutions are extinct at $t = 1$ under their flow), but this is violating the maximum principle which implies that the distance of two disjoint submanifold, if one of them is compact (which $M$ is), is non-decreasing under the flow - this gives a contradiction. Thus the boundary of $\bar N_{ext}$ is connected and the statement is true. 
\end{pf}

Of course, in the one point compactification, since $M$ is compact we see as sets $\widetilde{N}_{int} = N_{int}$. $N_{ext}$ is related to $\widetilde{N}_{ext}$ by way of the following observation:
\begin{lem} Suppose that $N \subset \R^3$ is a 3 manifold with boundary that is the exterior of a closed compact surface $M$, and let $\widetilde{N}$ be the compactification of $N$ induced by the one point compactification of $\R^3$ by adding a point. Then the induced map on $\pi(N) \to \pi(\widetilde{N})$ by inclusion is surjective. 
\end{lem}
\begin{pf} This is because $N$ contains a neighborhood of spatial infinity because $K$ is compact, and because the 2-sphere is simply connected.
\end{pf}
$\medskip$

Thus, arguing as in the proof of Theorem 3 in \cite{L} one can show that $M$ is a Heegaard surface and so is isotopic to the standard surface of genus $g$ in $S^3$ by Waldhausen's theorem \cite{W}, so any two compact self shrinkers $M_0$, $M_1$ of genus $g$ in $\R^3$ are ambiently isotopic to the standard surface of genus $g$ in the one point compactification of $\R^3$. 
$\medskip$

In addition, we may arrange this isotopy to avoid $\{\infty\}$ (the north pole) in $S^3$. To see this, consider an isotopy $M_t$ of $M_0$ to $M_1$. Since $M_0$ and $M_1$ are compact there is a $r > 0$ so that $M_0$ and $M_1$ are both disjoint from $B(\{\infty\}, r)$ as subsets of $S^3$. Note in fact (identifying . For each fixed $M_t$, we may arrange a rigid $R_t$, $R_0 = Id$, so that $R_t \circ M_t$ avoids $B(\{\infty\}, r)$. Since $M_t$ vary continuously the $R_t$ can be taken to as well, so we get a new isotopy $M_0$ to $M_2$, where $M_1 = R_1^{-1} M_2$, which avoids $B(\{\infty\}, r)$. Since $M_1 \cap B(\{\infty\}, r) = \emptyset$ as well we may may arrange an isotopy through rotations from $M_2$ to $M_1$ which avoids $B(\{\infty\}, r)$ (this is just since $S^3$ minus a ball is diffeomorphic to $\R^3$). 
$\medskip$

Hence by stereographic projection from $\{\infty\}$, the surfaces are isotopic in $\R^3$, giving theorem \ref{main}.

\section{Concluding Remarks.}
As hinted at above, analogous results have been found by Lawson in \cite{L} for minimal hypersurfaces in $S^3$ endowed with the round metric (to prove essentially lemma 2.1 above he used a Frenkel argument \cite{Fr} crucially using positive Ricci curvature) and later in general 3 manifolds manifolds with positive Ricci curvature by Meeks in \cite{M}. There are also results by Freedman, Frohman, Meeks, Yau, and others for different classes of minimal hypersurfaces in $\R^3$ - see \cite{MY}, \cite{FM}, \cite{F}. 
$\medskip$

In the conclusion of his paper Lawson gives (or at least restates) his famous conjecture (solved finally by Brendle in \cite{B}) that the Clifford torus is the only embedded minimal torus in $S^3$, his result (in particular, that there are no minimal knotted tori in $S^3$) meant as a step towards this. Similarly, we feel inclined to (re)state the following natural conjecture:
\begin{conj} The only embedded self shrinking torus in $\R^3$ is Angenent's torus. 
\end{conj} 
This conjecture, which has been stated before in various formulations (see for example the remark after theorem 2.2.2 of Moller's thesis \cite{Mo}) is not completely unfounded. Angenent used a ``shooting method" to construct a curve that, when rotated about an axis, gives a self-shrinker. This method suggests at least local isolatedness but the authors presently don't know how to prove this. We might also note (to be taken perhaps as evidence to the contrary) that Angenent in \cite{Ang}, the same paper where he constructs his famous embedded torus, constructs many immersed (not embedded) tori and Drugan and Kleene in \cite{DK} also constructed many immersed tori (they also have a construction of mebedded tori, and it is not clear if their tori is the same as Angenent's even). Of course, there are many immersed CMC surfaces but only the spheres are embedded, so perhaps this evidence isn't so threatening. 
$\medskip$

We end with noting that, while it is easy to see the Gaussian metric on $\R^n$ is not Ricci positive, it is $f$-Ricci positive for $f=\frac{|x|^2}{4}$, up to scaling, in the sense of Bakry-Emery. Using totally different techniques from the above, Wei and Wylie show that any two compact $f$-minimal hypersurfaces in a manifold with $Ric_f > 0$ intersect (more precisely, see theorem 7.4 in \cite{WW}) by showing the distance between them must be zero provided the distance between them is realized by the length of a smooth geodesic between two points on the surfaces (indeed, if there are two points realizing the distance it's fine if the minimal surfaces involved are noncompact - in the setting of \cite{L} Lawson achieves this by deck transformations to bring points in a distance minimizing sequence back to a fixed ball, see observation 3 in his paper). We point out however that since the Gaussian metric degenerates at infinity, which is within $N_{ext}$ and is a finite distance away in the Gaussian metric, its not clear we can directly apply the argument of Lawson on the universal cover of $N_{ext}$, $\overline{N_{ext}}$, using the theorem of Wei and Wylie to get our $\pi_1$ surjectivity and hence unknottedness statement. 
$\medskip$

More precisely, it is unclear how the degeneracy at infinity of the metric affects how bad the (in)completeness of the lifted metric is - the lift of the a neighborhood of infinity could concievably be ``between" two connected components of $\partial \overline{N_{ext}}$ and thus be relevent. Similarly we can't directly apply their argument because of the incompleteness to the end of the proof of lemma 2.1, but because the incompleteness of the Gaussian metric is only at infinity and one of the self shrinkers involved is compact one can perhaps alternatively use the method of Wei and Wylie to finish the proof.

\end{document}